\documentclass{article}
\usepackage[T1]{fontenc}
\usepackage[latin1]{inputenc}
\usepackage{geometry}
\makeatletter

\providecommand{\LyX}{L\kern-.1667em\lower.25em\hbox{Y}\kern-.125emX\@}

\usepackage[T1]{fontenc}
\usepackage[latin1]{inputenc}

\makeatletter

\usepackage[T1]{fontenc}
\usepackage[latin1]{inputenc}
\usepackage{geometry}
\makeatother
\everyjob{\typeout{MuTeX.sty: LaTeX Modification by EPM.}}
\immediate\write10{MuTeX.sty: LaTeX Modification by EPM.}

\makeatletter   % enable @ as normal letter;
\input{ifthen.sty}

\newtheorem {theorem} {Theorem} [section]

\newtheorem {lemma} [theorem] {Lemma}

\newtheorem {Canham threshold} [theorem] {Canham Threshold}

\def\theoremstyle#1#2{\def\@@theoremheadstyle{#1}
                      \def\@@theorembodystyle{#2}}
\def\@@theoremheadstyle{\sc}
\def\@@theorembodystyle{\rm}
\def\@begintheorem#1#2{\@@theorembodystyle 
                       \trivlist 
		       \item[\hskip 
                             \labelsep{\@@theoremheadstyle #1\ #2}]}
\def\@opargbegintheorem#1#2#3{\@@theorembodystyle 
                              \trivlist 
			       \item[\hskip 
				  \labelsep{\@@theoremheadstyle #1\ #2\ (#3)}]}

 \def\@@pc{\bf}

 \newcommand {\pcodestyle}[1] {\def\@@pc{#1}}  % use this to change the style

 {
 \def\PROGRAM		{{\@@pc program\ }}
 
 \def\PROCEDURE		{{\@@pc procedure\ }}
 \def\FUNCTION		{{\@@pc function\ }}
 \def\LOCAL		{{\@@pc local\ }}
 \def\GLOBAL		{{\@@pc global\ }}
 \def\RETURNS		{{\@@pc returns\ }}
 \def\RETURN		{{\@@pc return\ }}
 \def\BEGIN		{{\@@pc begin\ }}
 \def\END		{{\@@pc end\ }}
 \def\IF			{{\@@pc if\ }}
 \def\THEN		{{\@@pc then\ }}
 \def\ELSE		{{\@@pc else\ }}
 \def\REPEAT		{{\@@pc repeat\ }}
 \def\UNTIL		{{\@@pc until\ }}
 \def\WHILE		{{\@@pc while\ }}
 \def\DO			{{\@@pc do\ }}
 \def\FOR		{{\@@pc for\ }}
 \def\TO			{{\@@pc to\ }}
 \def\DOWN		{{\@@pc down\ }}
 \def\NEXT		{{\@@pc next\ }}
 \begin{center}
  \begin{minipage}{.9\textwidth}
   \small
    \begin{tabbing}
 }%
 {  \end{tabbing}%
   \end{minipage}%
  \end{center}%
 }

			   {\end{quotation}}

			 {  \end{minipage}%
			   \end{center}%
			  \end{description}%
			 }

				{ \end{minipage}
				 \end{center}%
				}

\def\thebibliography#1{\section*{References}\list
 {[\arabic{enumi}]}{\settowidth\labelwidth{[#1]}\leftmargin\labelwidth
 \advance\leftmargin\labelsep
 \usecounter{enumi}}
 \def\newblock{\hskip .11em plus .33em minus -.07em}
 \sloppy
 \sfcode`\.=1000\relax}

\newsavebox{\ProofSym}
\savebox{\ProofSym}{%
  \begin{picture}(10,10)
    \put(0,0){\framebox(9,9){}}
    \put(0,3){\framebox(6,6){}}
  \end{picture}}
\newcommand{\eop}{\hfill\usebox{\ProofSym}}
\newenvironment{proof}{\noindent {\sc Proof.\/}}{\eop\par\vspace{0.3cm}}

\begin{document}

\title{Computing sharp and scalable 
bounds on errors in approximate zeros of
univariate polynomials.
\footnote{This work was partially funded by a research grant 
from the AICTE, India.}}

\author{P. H. D. Ramakrishna\footnote{Currently, 
Deputy Superintendent of Police (Probationer), 
Andhra Pradesh, India}, 
Sudebkumar Prasant Pal\footnote{Department of Computer 
Science and Engineering,
Indian Institute of Technology, Kharagpur, 721302, India.
email: spp@cse.iitkgp.ernet.in. 
http://www.angelfire.com/or/sudebkumar} \footnote{A part of this work was 
done when this author was visiting
the University of Miami, Coral Gables, Florida, USA.}, 
Samir Bhalla,\\ 
Hironmay Basu, 
Sudhir Kumar Singh 
}
\date{}
\maketitle

\begin{abstract}

There are several numerical methods for computing approximate zeros 
of a given univariate polynomial. 
In this paper, we develop a simple and novel 
method for determining sharp upper 
bounds on errors in approximate zeros of a given polynomial 
using Rouche's theorem from complex analysis.
We compute the error bounds using non-linear optimization. Our bounds
are scalable in the sense that we compute sharper error bounds for
better approximations of zeros. 
We use high precision computations using the LEDA/real
floating-point filter 
for computing our bounds robustly.

\end{abstract}


\section{Introduction}

The problem of solving the polynomial equation 
\begin{equation}
\label{1}
g(z)=a_{0}+a_{1}z+a_{2}z^{2}+...+a_{n}z^{n}=0
\end{equation}
is a fundamental problem in algebra. Here, the coefficients 
$a_i$, $0\leq i\leq n$
are complex numbers. 
It is now folklore that no closed form formula with arithmetic 
operations and radicals
exists for computing zeros of 
polynomials with degree \( n>4. \) 
In spite of the absence of solution formulae, the \( fundamental \)
\( theorem \) \( of \) \( algebra \) states that equation \ref{1} always
has \( n \) zeros, real or complex. The subject of computing or approximating
these zeros has been called \emph{\( algebraic \) \( aspects \) \( of \)
\( the \) \( fundamental \) \( theorem \) \( of \) 
\( algebra \)} \cite{smale}.
The computational challenge is in determining zeros precisely for high degree
univariate polynomials.  
In computer algebra applications, one usually 
needs to solve equation \ref{1} for
large values of \( n \), typically well 
above \( 100 \) or of order of several thousands \cite{pan2}.
High precision of hundreds of bits is frequently required for the
representation of the coefficients \( a_{0},a_{1},...,a_{n} \) 
and the zeros.
Development of robust and efficient  solutions for
equation \ref{1} with large $n$ 
requires much research.

\subsection{Computing approximate zeros of polynomials}

Since there is no closed form formula for finding zeros 
of a polynomial of degree greater
than four, we need to find zeros by numerical methods. 
The first iterative algorithms with guaranteed convergence to all the \( n \)
zeros of \( g(z) \) (for any input polynomial $g(z)$ 
of degree \( n \) ), are due
to Brouwer and  Weyl, both 
published in 1924 (see \cite{pan2}). 
Pan's new \( O(n^{2}\log n\log bn) \) algorithm to approximate
all the zeros of a $n$-degree polynomial within 
an error bound of \( 2^{-b}max_{j}|z_{j}|, \)
where \( z_{j} \) are the zeros of the polynomial \cite{pan1}, is
the best known result so far. 
There are 
several algorithms for finding zeros of polynomials based on 
Newton's method \cite{M73},
Laguerre's method \cite{smith0}, Jenkins and Traub's method \cite{JT72}, 
and Muller's method \cite{pan2}.

\subsection{Estimating upper bounds on errors in approximate zeros 
of a given polynomial}

Since all the above methods essentially compute 
only approximate zeros, we  
need to know how close these approximations are to
the actual zeros. Consider a practical scenario where we have pretty good
estimates of all the zeros of a known univariate polynomial as in equation 
\ref{1}. Although the zeros may be quite accurate, we still need to 
know how good they are. In other words, we need to know how much error
is present in each approximate zero.
Consider the problem of developing a 
method for determining upper bounds on the errors
in each approximate zero,  
irrespective of the method and computational precision
by which such approximate zeros have been determined apriori.
Smith \cite{smith} computes upper bounds on errors in approximate
zeros using a method based on Gerschgorin's theorems.
In this paper we develop our new and novel method for finding
such posteriori 
upper bounds on errors in given approximate zeros of a given univariate
polynomial. 
Our method computes an upper bound on 
the {\it distance} of an approximate
zero from the exact zero
given all the coefficients of the univariate polynomial and  
all the approximate zeros.
Note that we are given only the 
univariate polynomial and its approximate zeros as inputs; we may not know 
any actual zero. Towards this end, we use 
Rouche's theorem from \( complex \) \( analysis \) \cite{ahlfors} as
the basis for our approximation scheme. 
Based on this scheme, we develop a non-linear optimization step to determine
a sharp upper bound on the modulus of the distance of an approximate zero 
from the (possibly unknown) exact zero. This method is reported in 
\cite{rama,techrep}. The implementation of our method for
estimating these error bounds uses the {\it LEDA/real floating-point filter} 
\cite{exactleda,mehlhornnaher99} so that computations can be done at requisite high 
precision
beyond double precision, yielding correct results in the evaluation of 
inequalities involving arithmetic expressions. The number of significant digits
in the approximate zeros and the coefficients of the given polynomial
are not restricted therefore.

\subsection{Preliminaries}

We need
the definition of an {\it analytic function} \cite{ahlfors}.
A {\it complex function} \( f(z) \) is said to be {\it analytic} 
in a region 
if it is defined and differentiable at every point $z$ in the region.
Now we state Rouche's theorem.
\begin{theorem}
\cite{ahlfors,henricivol1}
Suppose the functions \( f(z) \) and \( g(z) \) are analytic 
inside and on a simple
closed curve $C$.  If $f$ and $g$ have no zeros on 
\( C \) and \( |f(z)-g(z)|<|f(z)| \) for all \( z \) on
\( C \),
then the functions \( f(z) \) and \( g(z) \) have the same number of
zeros inside \( C \).
\label{rouchetheorem}
\end{theorem}
We know that all polynomials are analytic in the complex plane. 
We know from the {\it fundamental theorem of algebra} that any 
polynomial with complex coefficients of degree \( n \) has \( n \) 
(possibly) complex zeros. Given a polynomial
of order \( n \) and its \( n \) approximate zeros, 
we wish to bound the error in each of the \( n \) aproximate zeros. 
We present our method in the next section.

\section{A novel 
posteriori error bound for an approximate zero}

\label{novelmethod}

Let the given polynomial be 
\begin{equation}
\label{2}
g(z)=a_{0}+a_{1}z+a_{2}z^{2}+...+a_{n}z^{n}
\end{equation}
where the coefficeints $a_i$, $0\leq i\leq n$ are complex in 
general.
Let the approximate zeros of \( g(z) \) (say, as computed by 
some numerical method) be \( \alpha _{1},\alpha _{2},\cdots ,\alpha _{n} \). 
Let 
\begin{equation}
\label{3}
f(z)=a_{n}\prod_{i=1}^n (z-\alpha _{i})
\end{equation}
 be the approximation of $g(z)$ as constructed from the given
approximate zeros.
Now consider the {\it error} polynomial \( h(z) =f(z)-g(z) \). 
Let 
\begin{equation}
\label{4}
h(z)=b_{0}+b_{1}z+.....+b_{n-1}z^{n-1}
\end{equation}

Note that the error polynomial $h(z)$ has degree at most $n-1$, less than
that of $g(z)$.
Now we apply Rouche's theorem on the polynomials \( f(z) \) and \( g(z). \) 
Consider \( \alpha _{j} \), an approximate zero. In this section and 
in section \ref{computingerror}, 
we concentrate on $\alpha_j$; the arguments in 
these sections apply to each approximate zero.
We draw a circle \( C \)
with a small radius \( r \) around \( \alpha _{j} \). If the error polynomial
\( h(z)=f(z)-g(z) \) and the computed polynomial \( f(z) \) satisfy 
Rouche's condition
\( |f(z)-g(z)|<|f(z)| \) for all \( z \) on the 
on this circle, then \( f(z) \) and \( g(z)  \)
have precisely the same number of zeros inside this circle. 
If there is no other
approximate zero \( \alpha _{i}, i\neq j \) within that circle, 
then \( f(z) \) and \( g(z) \)
have precisely one zero in the circle. So, the unique actual zero of $g(z)$ is
approximated by \( \alpha _{j} \) and is within that circle of small radius 
$r$ with centre at $\alpha_j$. Here 
$r$ is an upper bound on the error in the approximate zero \( \alpha _{j} \).

For enforcing the condition in Rouche's theorem we observe that,
for any $z$ on the circle $C$ of radius $r$,
\( |z|=|z-\alpha _{j}+\alpha _{j}|\leq |z-\alpha _{j}|+|\alpha _{j}| \) = \( r+|\alpha _{j}| \)
= \( R \), say.
From equation \ref{4}, \( |h(z)|\leq |b_{n-1}||z^{n-1}|\cdots +|b_{0}| \). 
Hence on \( C, \) 
\begin{equation}
\label{5}
|h(z)|\leq |b_{n-1}| R^{n-1}+\cdots +|b_{0}|
\end{equation}
Also, for all $z$ on the circle \( C \) of radius $r$,
\( |z-\alpha _{i}|=|(z-\alpha _{j})-(\alpha _{i}-\alpha _{j})|\geq 
||z-\alpha _{j}|-|\alpha _{i}-\alpha _{j}|| \)
= \(| r-|\alpha _{i}-\alpha _{j}|| \).
Also, \( |f(z)|=|a_{n}||z-\alpha _{1}||z-\alpha _{2}|\cdots |z-\alpha _{n}| \),
from equation \ref{3}.
Hence on $C$, 
\begin{equation}
\label{6}
|f(z)|\geq |a_{n}|r\prod_{i=1, i\neq j}^{n} |r-|\alpha _{i}-\alpha _{j}||
\end{equation}

Finally, from inequalities \ref{6} and \ref{5},
we observe that  \( |a_{n}|r|\prod_{i=1, i\neq j}^{n} |r-|\alpha _{i}-\alpha _{j}||>|b_{n-1}|R^{n-1}+.....+|b_{0}| \) implies
Rouche's condition \( |f(z)-g(z)|<|f(z)| \) for all $z$ 
on the entire circle $C$.
We can write this condition as
\begin{equation}
\label{7}
r>\frac{l(r)}{m(r)}(=q(r), say)
\end{equation}
where $ l(r)= |b_{n-1}|R^{n-1}+....+|b_{0}|$ and $ m(r)= |a_{n}|\prod_{i=1, i\neq j}^n |r-|\alpha _{i}-\alpha _{j}||$.
Note that $l(r)$ and $m(r)$ are of same degree $n-1$ and  
that $q(r)$ is a {\it meromorphic} function.

Now we summarize the main theorem of this paper as follows.
\begin{theorem}
\label{maintheorem}
Let $g(z)$ be a polynomial of degree $n$ with $n$ distinct zeros. Let
\( \alpha _{1},\alpha _{2},\cdots ,\alpha _{n} \) be approximations
to the $n$ zeros of $g(z)$ and let 
$f(z)=a_{n}\prod_{i=1}^n (z-\alpha _{i})$
be the polynomial with zeros
\( \alpha _{1},\alpha _{2},\cdots ,\alpha _{n} \). Then, for
any $1\leq i\leq n$, the error in the $j$th approximate zero $\alpha_j$ is
bounded by any real number $r$ satisfying the inequality $r>q(r)$, provided
$|\alpha_j - \alpha_i| > r$, for all $1\leq i\leq n$ where $i\neq j$.
\end{theorem}
\begin{proof}
As shown above, 
Rouche's condition is satisfied on the circle $C$ centred at $\alpha_j$ with 
radius $r$ provided $r>q(r)$. If no other approximate zero $\alpha_i$, 
$i\neq j$ is in the interior of 
this circle $C$ then $f(z)$ has a single zero 
$\alpha_j$ in the interior of $C$. By Theorem \ref{rouchetheorem}, 
$g(z)$ too has a
single zero (say, $z_j$) inside the circle $C$, 
yielding the upper bound $r$ on the error
$|\alpha_j     -z_j|$.
\end{proof}

\section{A method for computing error bounds}

\label{computingerror}

We state the inequality \ref{7} as 
\begin{equation}
\label{8}
\frac{l(r)}{r*m(r)}<1
\end{equation}

We observe that the  inequality  \ref{8}  (and therefore  the inequality 7) 
is indeed satisfied for sufficiently large values of \(r \); 
this is due to the fact that \(r*m(r) \) 
is a degree \(n\) polynomial and \(f(r) \) is degree 
\( n-1 \) polynomial. 
We state this a follows.
\begin{lemma}
There exists a value or $r$ satisfying inequality \ref{8}.
\label{largebound}
\end{lemma}
However, we are interested in finding very 
small values of \( r\) satisfying \( r > q(r)\). In particular,
note that in our Theorem \ref{maintheorem}, such an \(r \)  is 
an upper bound on the error in an approximate zero $\alpha_j $
provided $|\alpha_j-\alpha_i| >r$, for all $i\neq j $. 
We proceed to develop a method to find such a value \( r\) as follows.

\subsection{Computing the error bounds: Algorithm I}
\label{algoI}

We argue that inequality \ref{7} 
is satisfied for values of $r$ close to $q(0)$, 
provided $f(z)$ is a close approximation of $g(z)$ and the given 
approximate zeros are well separated. 
If $f(z)$ is close to $g(z)$, then $l(0)$ is small and therefore $q(0)$ could 
be small.  By Taylor's 
expansion, \( q(r)=q(0)+rq'(0)+\frac{r^{2}}{2!}q''(0)+\cdots  \).
Since higher derivatives of $q(r)$ will be small for well separated 
approximate zeros, we can neglect 
higher order terms and write the fixpoint of the function $q$ as
\( r=q(r)=q(0)+rq'(0) \) 
\( \Rightarrow r=\frac{q(0)}{1-q'(0)}. \)
Note that \( q'(0) \) is also small since \( q(r) \) is a function
with a numerator $l(r)$ with very small coefficients, and 
a denominator $m(r)$ with 
large coefficients. Therefore, the fixpoint of $q$ lies in the
vicinity of \( q(0) \), provided the approximate zeros are well 
separated and $f(z)$ is a good approximation of $g(z)$.

{\tiny

\centering
\subsection*{Example 1}

\centerline{\( g(z)=100000z^{4}+305000z^{3}+410100z^{2}+310205z+105105 \)
}

\par{} \vspace{0.3cm}\centering
Comparision  with Smith's bounds

\vspace{0.3cm}\centering
\begin{tabular}{|c|c|c|c|c|}
\hline 
Actual &
 Zeros by&
 Smith's bounds&
 Value &
 Our bounds on \\
 zeros&
 ZERPOL&
on ZERPOL &
 of \( q(0) \)&
 ZERPOL zeros\\
&
&
 zeros&
&
 using Algorithm I\\
&
&
&
&
\( \epsilon =0.00001 \)\\
\hline 
\( -1.05 \)&
 \( -1.0500001610 \)&
 \( 6.44E-07 \)&
 \( 1.116840490 \)&
 \( 1.11687399 \)\\
&
&
&
 \( 680863E-06 \)&
 \( 61584E-06 \)\\
\hline 
\( -1 \)&
 \( -0.9999998510 \)&
 \( 5.97E-07 \)&
 \( 1.08266830 \)&
 \( 1.08270078 \)\\
&
&
&
 \( 781336E-06 \)&
 \( 81847E-06 \)\\
\hline 
\( -0.5+i\sqrt{0.751} \)&
\multicolumn{1}{c|}{ \( -0.5+i0.86660 \)}&
 \( 2.55E-08 \)&
 \( 3.04704508 \)&
 \( 3.04707641 \)\\
&
\multicolumn{1}{c|}{\( 2562368 \)}&
&
 \( 007204E-08 \)&
 \( 434807E-08 \)\\
\hline 
\( -0.5-i\sqrt{0.751} \)&
 \( -0.5-i0.86660 \)&
 \( 2.55E-08 \)&
 \( 3.04704594 \)&
 \( 3.04707641 \)\\
&
 \( 2562368 \)&
&
 \( 388862E-08 \)&
 \( 434806E-08 \) \\
\hline 
\end{tabular}

\par{} \vspace{0.7cm}\centering

\vspace{0.3cm}
{\centering \begin{tabular}{|c|}
\hline 
Zeros at precision 7\\
\hline 
\end{tabular}\par}
\vspace{0.3cm}

Bounds using Algorithm I for \( \epsilon  \)\( =0.0001 \)

\begin{tabular}{|c|c|c|c|}
\hline 
Range zeros at &
 Value of \( q(0) \)&
 Our bounds on Range&
 Number of \\
 precision 7&
&
 zeros at precision 7&
 iterations\\
\hline 
\( -1.05 \)&
 \( 5.42072490014779E-06 \)&
 \( 5.42180887902023E-06 \)&
 \( 2 \)\\
\hline 
\( -1.000000 \)&
 \( 5.43336060444733E-06 \)&
 \( 5.4344471 \)\( 1286507E-06 \)&
 \( 2 \)\\
\hline 
\( -0.5+i0.8666026 \)&
 \( 1.5286312353107E-07 \)&
 \( 1.52878409 \)\( 843423E-07 \)&
 \( 1 \)\\
\hline 
\( -0.5-i0.8666026 \)&
 \( 1.52863125747657E-07 \)&
 \( 1.52878412 \)\( 060232E-07 \)&
 \( 1 \) \\
\hline 
\end{tabular}

\par{} \vspace{0.7cm}\centering

Bounds using Algorithm II

\begin{tabular}{|c|c|c|c|}
\hline 
Range zeros at&
 Starting &
 NR bounds on Range&
\multicolumn{1}{c|}{Number of iterations}\\
 precision 7&
 Value&
 zeros at precision7&
\multicolumn{1}{c|}{\( NR \) \( + \)Algorithm I}\\
\hline 
\( -1.05 \)&
 \( 0.01 \)&
 \( 5.42194046573675E-06 \)&
 \( 4+1 \)\\
\hline 
\( -1.000000 \)&
 \( 0.001 \)&
 \( 5.43458277061677E-06 \)&
 \( 4+1 \)\\
\hline 
\( -0.5+i0.8666026 \)&
 \( 100 \)&
 \( 1.52878500883089E-07 \)&
 \( 4+1 \)\\
\hline 
\( -0.5-i0.8666026 \)&
 \( 1E-05 \)&
 \( 1.5287850088089E-07 \)&
 \( 3+1 \) \\
\hline 
\end{tabular}

\par{} \vspace{0.7cm}

\vspace{0.3cm}
{\centering \begin{tabular}{|c|}
\hline 
Zeros at precision 16\\
\hline 
\end{tabular}\par}
\vspace{0.3cm}

Bounds using Algorithm I for \( \epsilon  \)\( =0.00000001 \)

\begin{tabular}{|c|c|c|c|}
\hline 
Range zeros at&
 Value of \( q(0) \)&
 Our bounds on Range&
 Number of\\
 precision 16&
&
 zeros at precision 16&
iterations\\
\hline 
\( -1.05 \)&
 \( 4.64160094696633E-15 \)&
 \( 4.64160099338234E-15 \)&
 \( 1 \)\\
\hline 
\( -1.000000000000000 \)&
 \( 4.51503828651945E-15 \)&
 \( 4.51503833166985E-15 \)&
 \( 1 \)\\
\hline 
\( -0.5+i0.866602561731732 \)&
 \( 1.27065229457289E-16 \)&
 \( 1.27065230727944E-16 \)&
 \( 1 \)\\
\hline 
\( -0.5-i0.866602561731732 \)&
 \( 1.27065229457292E-16 \)&
 \( 1.27065230727944E-16 \)&
 \( 1 \) \\
\hline 
\end{tabular}

\par{} \vspace{0.7cm}\centering

Bounds using Algorithm II

\begin{tabular}{|c|c|c|c|}
\hline 
Range zeros at&
 Starting &
 NR bounds on Range&
 Number of iterations\\
 precision 16&
 Value &
 zeros at precision 16&
\multicolumn{1}{c|}{\( NR \) \( + \)Algorithm I}\\
\hline 
\( -1.05 \)&
 \( 1E-10 \)&
 \( 4.64160099338286E-15 \)&
 \( 2+1 \)\\
\hline 
\( -1.000000000000000 \)&
 \( 0.0001 \)&
 \( 4.51503833167034E-15 \)&
 \( 3+1 \)\\
\hline 
\( -0.5+i0.866602561731732 \)&
 \( 1E-22 \)&
 \( 1.27065230727944E-16 \)&
 \( 2+1 \)\\
\hline 
\( -0.5-i0.866602561731732 \)&
 \( 1E-14 \)&
\( 1.27065230727944E-16 \) &
 \( 2+1 \) \\
\hline 
\end{tabular}

\par{} \vspace{0.7cm}\centering

\subsection*{Example 2}

\vspace{0.3cm}

\( g(z)=1000z^{10}-2500z^{9}-460800z^{8}-9133400z^{7}-50761800z^{6}-88653100z^{5}- \)

\( 53510400z^{4}-37313000z^{3}-197170000z^{2}-364800000z-198000000 \)

\par\vspace{0.3cm}\centering
Comparision with Smith's bounds
\par\vspace{0.3cm}\centering
\begin{tabular}{|c|c|c|c|c|}
\hline 
Actual zeros&
 Zeros by &
 Smith's bounds &
Value of&
Our bounds on \\
&
ZERPOL&
on ZERPOL &
\( q(0) \)&
ZERPOL zeros\\
&
&
zeros&
&
using Algorithm I\\
&
&
&
&
\( \epsilon =0.00001 \)\\
\hline 
\( 30 \)&
\( 30 \)&
\( 3.34E-13 \)&
\( 6.16818594 \)&
\( 6.1707644 \)\\
&
&
&
\( 879537E-11 \)&
\( 473737E-11 \)\\
\hline 
\( -10+i10 \)&
 \( -10+i10 \)&
 \( 8.3E-13 \)&
\( 3.99156437 \)&
 \( 3.99160428 \)\\
&
&
&
\( 089776E-10 \)&
\( 654147E-10 \)\\
\hline 
\( -10-i10 \)&
 \( -10-i10 \)&
 \( 8.3E-13 \)&
\( 3.99156437 \)&
\( 3.9916042 \) \\
&
&
&
\( 089776E-10 \)&
\( 8654147E-10 \)\\
\hline 
\( -5 \)&
 \( -5 \)&
 \( 2.91E-10 \)&
\( 3.18914558 \)&
 \( 3.18914637 \)\\
&
&
&
\( 60219E-09 \)&
\( 774777E-09 \)\\
\hline 
\( 1+i \)&
 \( 1+i \)&
 \( 6.81E-14 \)&
\( 3.72504031 \)&
 \( 3.7250775 \)\\
&
&
&
\( 303826E-10 \)&
\( 6344139E-10 \)\\
\hline 
\( 1-i \)&
 \( 1-i \)&
 \( 6.81 \)\( E-14 \)&
\( 3.72504031 \)&
\( 3.7250775 \)\\
&
&
&
\( 303826E-10 \)&
\( 6344139E-10 \)\\
\hline 
\( -1+i\sqrt{1.2} \)&
 \( -1+i1.0954 \)&
 \( 3.62E-08 \)&
\( 2.9138084 \)&
 \( 2.91383760 \)\\
&
\( 451114 \)&
&
\( 647204E-09 \)&
\( 280512E-09 \)\\
\hline 
\( -1-i\sqrt{1.2} \)&
 \( -1-i1.0954 \)&
 \( 3.62E-08 \)&
\( 2.9138084 \)&
\( 2.9138376 \) \\
&
\( 451114 \)&
&
\( 647204E-09 \)&
\( 0280512E-09 \)\\
\hline 
\( -1.5 \)&
 \( -1.5 \)&
 \( 1.68E-13 \)&
\( 1.17048091 \)&
\( 1.1704926 \)\\
&
&
&
\( 956553E-08 \)&
\( 2437473E-08 \)\\
\hline 
\( -1 \)&
 \( -1 \)&
 \( 8.25E-14 \)&
\( 8.1388826 \)&
 \( 8.1389640 \)\\
&
&
&
\( 9390312E-09 \)&
\( 8273006E-09 \)\\
\hline 
\end{tabular}

\par{} \vspace{0.3cm}

\pagebreak

{\centering \begin{tabular}{|c|}
\hline
Zeros at precision 7\\
\hline
\end{tabular}\par}

Bounds from Algorithm I for \( \epsilon =0.0001 \)

\begin{tabular}{|c|c|c|c|}
\hline 
Range zeros&
Value of &
Our bounds on Range&
Number of\\
at precision 7&
\( q(0) \)&
zeros at precision 7 &
iterations\\
\hline 
\( 30 \)&
\( 7.02434379707557E-09 \)&
\( 7.02443898949365E-09 \)&
\( 1 \)\\
\hline 
\( -10+i10 \)&
\( 4.54376449207249E-08 \)&
 \( 4.54380992971741E-08 \)&
\( 1 \)\\
\hline 
\( -10-i10 \)&
\( 4.54376449207249E-08 \)&
\( 4.54380992971741E-08 \) &
\( 1 \)\\
\hline 
\( -5 \)&
\( 3.62210633642892E-07 \)&
 \( 3.62214255749229E-07 \)&
\( 1 \)\\
\hline 
\( 1+i \)&
\( 4.24036915436498E-08 \)&
 \( 4.24041155805652E-08 \)&
\( 1 \)\\
\hline 
\( 1-i \)&
\( 4.24036915436498E-08 \)&
\( 4.24041155805652E-08 \)&
\( 1 \)\\
\hline 
\( -1+i1.095445 \)&
\( 3.31691055907562E-07 \)&
 \( 3.31694372818121E-07 \)&
\( 1 \)\\
\hline 
\( -1-i1.095445 \)&
\( 3.31691055907562E-07 \)&
\( 3.31694372818121E-07 \) &
\( 1 \)\\
\hline 
\( -1.5 \)&
\( 1.33240756345005E-06 \)&
 \( 1.33242088752569E-06 \)&
\( 1 \)\\
\hline 
\( -1 \)&
\( 9.26483214486361E-07 \)&
 \( 9.26492479318505-07 \)&
\( 1 \)\\
\hline 
\end{tabular}

\par{} \vspace{0.3cm}

\par{} \vspace{0.3cm}

Bounds from Algorithm II

\begin{tabular}{|c|c|c|c|}
\hline 
Range zeros at&
Starting &
NR bounds on Range&
Number of\\
 precision 7&
value&
zeros at precsion 7&
iterations\\
&
&
&
\( NR \) +Algorithm I\\
\hline 
\( 30 \)&
\( 1E-5 \)&
 \( 7.02507120754467E-09 \)&
\(3+1 \)\\
\hline 
\( -10+i10 \)&
\( 1E-1 \)&
\( 4.54421907500478E-08 \)&
\( 4+1 \)\\
\hline 
\( -10-i10 \)&
\( 1E-10 \)&
 \( 4.54421907500478E-08 \)&
\(3+1 \)\\
\hline 
\( -5 \)&
\( 1E-1 \)&
  \( 3.62247181860654E-07 \)&
\( 4+1 \)\\
\hline 
\( 1+i \)&
\( 1E-5 \)&
 \( 4.24079393150912E-08 \)&
\( 3+1 \)\\
\hline 
\( 1-i \)&
\( 1E-5 \)&
\( 4.24079393150912E-08 \)&
\( 3+1 \)\\
\hline 
\( -1+i1.095445 \)&
\( 1E-10 \)&
 \( 3.317247877106E-07 \)&
\( 3+1 \)\\
\hline 
\( -1-i1.095445 \)&
\( 1E-2 \)&
 \( 3.317247877106E-07 \) &
\( 3+1 \)\\
\hline 
\( -1.5 \)&
\( 1E-2 \)&
 \(  1.33255239212745E-06 \)&
\( 4+1 \)\\
\hline 
\( -1 \)&
\( 1E-2 \)&
  \( 9.26581422929994E-07 \)&
\( 4+1 \)\\
\hline 
\end{tabular}
\par{} \vspace{0.3cm}

{\centering \begin{tabular}{|c|}
\hline
Zeros at precision 16\\
\hline
\end{tabular}\par}

Bounds from Algorithm I for \( \epsilon =0.00000001 \)

\begin{tabular}{|c|c|c|c|}
\hline 
Range zeros at&
 Value of&
Our bounds on Range &
Number of\\
precision 16&
\( q(0) \)&
zeros at precision 16&
iterations\\
\hline 
\( 30 \)&
\( 1.83299375288556E-17 \)&
 \( 1.83299377121549E-17 \)&
\( 1 \)\\
\hline 
\( -10+i10 \)&
\( 1.18568546630189E-16 \) &
 \( 1.18568547815876E-16 \)&
\( 1 \)\\
\hline 
\( -10-i10 \)&
\( 1.18568546630189E-16 \)&
\( 1.18568547815876E-16 \) &
\( 1 \)\\
\hline 
\( -5 \)&
\( 9.41580760246742E-16 \) &
 \( 9.45180769698553E-16 \)&
\( 1 \)\\
\hline 
\( 1+i \)&
\( 1.10651507753055E-16 \) &
\( 1.10651508859572E-16 \) &
\( 1 \)\\
\hline 
\( 1-i \)&
\( 1.10651507753055E-16 \) &
\( 1.10651508859572E-16 \) &
\( 1 \)\\
\hline 
\( -1+i1.095445 \)&
\( 8.65540430075807E-16 \) &
\( 8.65540438731223E-16 \)&
\( 1 \)\\
\( 115010332 \)&
&
&
\\
\hline 
\( -1-i1.095445 \)&
\( 8.65540430075807E-16 \) &
\( 8.65540438731223E-16 \)&
\( 1 \)\\
\( 115010332 \)&
&
&
\\
\hline 
\( -1.5 \)&
\( 3.47688796630911E-15 \)&
 \( 3.47688800107798E-15 \)&
\( 1 \)\\
\hline 
\( -1 \)&
\( 2.41763729896059E-15 \) &
 \( 2.41763732313696E-15 \) &
\( 1 \)\\
\hline 
\end{tabular}

\par{} \vspace{0.3cm}

Bounds from Algorithm II

\begin{tabular}{|c|c|c|c|}
\hline 
Range zeros at&
 Starting&
NR bounds on Range &
Number of \\
precision 16&
value&
zeros at precision 16&
iterations\\
&
&
&
\( NR \) +Algorithm I\\
\hline 
\( 30 \)&
\( 1E-14 \)&
\( 1.83299377120949E-17 \)&
\( 2+1 \)\\
\hline 
\( -10+i10 \)&
 \( 1000 \)&
\( 1.1856854781587E-16 \) &
\( 3+1 \)\\
\hline 
\( -10-i10 \)&
\( 1E-04 \)&
\( 1.1856854781587E-16 \)&
\( 3+1 \)\\
\hline 
\( -5 \)&
 \( 1E-04 \)&
\( 9.45180769698555E-16 \)&
\( 3+1 \)\\
\hline 
\( 1+i \)&
\( 1E-14 \) &
 \( 1.10651508859572E-16 \)&
\( 2+1 \)\\
\hline 
\( 1-i \)&
\( 1E-14 \) &
\( 1.10651508859572E-16 \) &
\( 2+1 \)\\
\hline 
\( -1+i1.095445115010332 \)&
\( 1E-14 \)&
 \( 8.65540438731227E-16 \)&
\(2+1 \)\\
\hline 
\( -1-i1.095445115010332 \)&
\( 1E-04 \) &
\( 8.65540438731227E-16 \)&
\( 3+1 \)\\ 
\hline 
\( -1.5 \)&
\( 1E-03 \) &
 \( 3.47688800107808E-15 \)&
\( 3+1 \)\\
\hline 
\( -1 \)&
\( 1E-03 \) &
\( 2.41763732313701E-15 \) &
\( 3+1 \)\\
\hline 
\end{tabular}

\par{} \vspace{0.3cm}

}

So, we start testing \( r \) with an initial value \( q(0). \) 
We increase \( r \) until \( r>q(r) \).
Note that $q(0)>0$,
and therefore, the inequality $r>q(r)$ is not satisfied at $r=0$. We
use $q(0)$ as a starting point for searching a value
of $r$ satisfying $r>q(r)$. In the search, we multiply the value of $r$
by a factor of $1+\epsilon$ repeatedly (for a small $\epsilon$), 
until the inequality $r>q(r)$ is satisfied. 
Further, if \( r<|\alpha _{j}-\alpha _{i}| \) 
for all \( 1\leq i\neq j\leq n \), we can assert using 
Theorem \ref{maintheorem} that $r$ is an upper bound on the error
in the approximate zero \( \alpha_j  \).

We claim that the number of times $q(r)$ is evaluated until we get a solution
$\delta$ satisfying $\delta > q(\delta)$ is $O(\log_2 \delta)$. This follows
from the fact that we start with an initial approximation $q(0)$ and multiply
$r$ by a factor $(1+ \epsilon)$ until the inequality \ref{7} is satisfied, as
guaranteed by Lemma \ref{largebound}.
We state this as a theorem. 
\begin{theorem}
The value of $r$ (say $\delta$), 
satisfying the inequality \ref{7} can be 
determined in $\log_{1+\epsilon} (\delta/q(0))$ evaluations of 
the function $q(r)$. 
\label{timecomplexity}
\end{theorem}
Each evaluation of $q(r)$ involves evaluation of 
two $(n-1)$-degree polynomials,
$l(r)$ and $m(r)$. Note in the tables (Tables 1 through 6)
that as long as we select 
a good starting value 
$q(0)$ for $r$, we practically need to execute just one step to get to a
value of $r$ satisfying $r>q(r)$.
Indeed, the error bounds are just a little bit higher than $q(0)$ in all cases,
very much as argued above. So, we observe that for good approximations
of well separated zeros (as considered in all our examples), 
Algorithm I converges very fast. Algorithm I requires larger number of steps 
if (i) the initial value of $r$ is far from a 
feasible solution for inequality 7 and,
(ii) for smaller values of $\epsilon$.

\subsection{Using requisite high precision for computing the error bounds}

Another factor in the time 
complexity of Algorithm I is the precision at which 
computations would require to be performed to correctly determine satisfaction
of inequality 7. We use the {\it LEDA/real floating-point filter}
\cite{exactleda,mehlhornnaher99}, for computations at requisite higher precisions.
LEDA first tries to check the 
inequality using double precision. If a decision can 
not be made at double precision, LEDA uses higher precision to check 
the inequality. Precision is increased until a decision is possible.
High precision computing is necessary
since round-off errors can accumulate in computations 
like the evaluation of high-degree polynomials. Such polynomials
may have coeffcients with
a large number of significant digits. Moreover, huge errors can accumulate
if a polynomial is evaluated at a
point whose value is a number with a large number of significant digits. 

\subsection{Computing the error bounds using the Newton-Raphson 
method: Algorithm II}

We can also solve inequality 7 using a combination of the Newton-Raphson 
(NR) method and Algorithm I. We call this method our Algorithm II. The main 
motivation here is to study the good behaviour of the function $p(r)=r-q(r)$ 
and solve the inequality quickly even if we use a large starting value for 
$r>>q(0)$. Our Algorithm II is as follows. 
Essentially, we use the function $p(r)=r-q(r)$ and solve for $p(r)=0$ using
NR iterations. We may start with a reasonably large value of $r$ so that 
$r>q(r)$ (or $p(r)>0$). We may also start with a smaller value of $r$ where 
$p(r)$ is negative. 
We perform NR iterations 
until the iterate $r$ does not change 
more than a small value say, $10^{-30}$.
The value of $p(r)$ for the terminating value 
of $r$ is observed to be negative in all our runs. 
Now we switch over to 
Algorithm I with this value of $r$ as starting value.
Note that Algorithm I always terminates when executed with such a
starting value of $r$ where $p(r)$ is negative. Finally 
at termination, a sharp small bound $r$ satisfying inequality 7 is obtained. 
The six tables show the results for six different polynomials. 

The function $p(r)$ responds well to the NR method. As mentioned earlier, 
$q'(r)$ is small (see Section \ref{algoI}), 
making $p'(r)$ close to 1 (note that $p'(r)=1-q'(r)$). 
Since $p'(r)$ is close to 1, NR has quadratic convergence.
We need only few NR iterations for each zero of each polynomial;
we find that 
the value of $r$ resulting after these 
small number of NR iterations is such 
that $p(r)$ is negative or equivalently, $r<q(r)$. As mentioned above, such a 
value of $r$ is a suitable starting point for Algorithm I.

Now we concentrate on the NR iterations in Algorithm II
required to approximately compute 
the fixpoint of $q$. Let $\rho=q(\rho)$ be the fixpoint.
Suppose we start with a large starting value $R$ for the NR iterate
$r$. Note that in a few steps we would quickly get a value of the 
iterate $r$ such that the first few significant bits of $r$ match
with the first few
significant bits of $\rho$. This happens because the first step corrects
the initial value $R$ to $R-p(R)/p'(R)$, 
which is nearly $q(R)$ and very small compared to $R$. This happens since
$p(r)=r-q(r)$, $p'(r)=1-q'(r)$ and $q'(R)$ is small eventhough $R$ may 
be quite large. So, in a few NR iterations $r$ gets close 
enough to $\rho$ to match $\rho$ in a few significant bits (binary digits).

Suppose $r$ and $\rho$ have $c\geq 1$ identical (most) significant bits. 
Due to quadratic convergence property of NR iterations, the number of 
matching significant bits will increase in geometric progression with a ratio
greater than unity. So, after $k$ iterations, $w^kc$ bits would match, where
$w>1$. If we are interested in getting $b$ bits of $r$ match $b$ bits of 
$\rho$, we would then need $k$ steps where $b=w^kc$. We assume without loss
of generality that $\rho=d\times 2^E$, where $0.5\leq d<1$ and $E$ is 
any integer. Then, in $k=O(\log_w b)$ iterations, we would get the iterate $r$
close to $\rho$ with relative error 
$|r-\rho|/\rho\leq (2^{-b}\times 2^E)/2^{E-1}=2^{-(b-1)}$.
We summarize this in the following theorem

\begin{theorem}
\label{nrtheorem}
The fixpoint $\rho=q(\rho)$ can be approximated within relative error 
$2^{-b}$ in $O(\log b)$ Newton-Raphson iterations, where $b>1$ is an integer.
\end{theorem}

\subsection{Coping with errors in the input process}

Another aspect that needs mention is the possibility of small 
errors creeping in
in the process of feeding inputs into our error computing algorithm. There 
are only two sets of inputs; the coefficients of the input polynomial and the
approximate zeros. We can avoid input errors in the coefficients if 
each coefficient is integral. If we have decimal rational numbers 
as coefficients, the input of such numbers may suffer an error in the 
process of standard input in any programming language. So,  
we multiply all the coefficients with a suitable large power of 10 to 
make all coefficients integral. This does not change the polynomial and 
therefore the zeros remain unchanged by such multiplication of the 
coefficients of the input polynomial. Indeed, for this very reason, 
all examples (except example 5) illustrated in this paper
are about polynomials with integral coefficients. 
The approximate zeros in the input 
may however originate from any source or any numerical method. 
In the context of this paper, we require that approximate zeros 
be propagated accurately into our 
error finding Algorithms I and II; 
if there is any error in the 
input of approximate zeros, these errors must be accounted for
in the final error bound computed for each zero. For the sake of simplicity,
we may assume that such errors do not 
occur in the input of the approximate zeros; 
after all, once the intended approximate zeros are 
entered into variables of our error computing program, the values 
of the approximate zeros and all coefficients of the input polynomial 
are processed considering all errors that 
may occur during the computation of the final error bounds on each approximate
zero. This is ensured by the floating-point filter and the {\it real}
data type in LEDA \cite{exactleda,mehlhornnaher99}.  
It is however possible and 
necessary to determine safe upper bounds in input processes for inputs from
standard input instructions in different programming languages and systems
and incorporate those error bounds suitably into the final error bound.
A detailed study of the issue of errors in input processes 
for computations of 
{\it safe} boolean operations between polygons 
is done in \cite{mpvt96}.  Naturally, computations of errors
resulting due to round-off errors in finite precision floating-point 
computations must also include due consideration of possible errors 
in input variables themselves; safe upper 
bounds on errors in each input variable must be provided and accounted for
in the propagation of errors through the entire computation until outputs are
generated. Such error analysis 
is reported in \cite{rama,koul}.   

{\tiny

\centering
\subsection*{Example 3}

\( g(z)=1000z^{6}-13016z^{5}+59214z^{4}-107974z^{3}+61769z^{2}-997z+4 \)

\vspace{0.3cm}\centering 

\vspace{0.3cm}
{\centering \begin{tabular}{|c|}
\hline 
Zeros at precision 7\\
\hline 
\end{tabular}\par}
\vspace{0.3cm}

Bounds using Algorithm I for \( \epsilon =0.0001 \)

\begin{tabular}{|c|c|c|c|}
\hline 
Range zeros &
Value Of &
 Our bounds on Range &
Number of \\
at precision 7&
\( q(0) \)&
zeros at precision 7&
iterations\\
\hline 
\( 4.993716 \)&
\( 2.089829937337557E-05 \)&
 \( 2.0899762298219E-05 \)&
\( 7 \)\\
\hline 
\( 4.0080655 \)&
\( 3.64330724956225E-05 \)&
 \( 3.64381734573261E-05 \)&
\( 14 \)\\
\hline 
\( 2.997721 \)&
\( 1.65965389611843E-05 \)&
 \( 1.65977007537648E-05 \)&
\( 7 \)\\
\hline 
\( 1.0 \)&
\( 7.32990974711188E-07 \)&
 \( 7.32998304620935E-07 \)&
\( 1 \)\\
\hline 
\( 0.009422683 \)&
\( 1.26751047717699E-08 \)&
 \( 1.26752315228176E-08 \)&
\( 1 \)\\
\hline 
\( 0.007075149 \)&
\( 9.44998978980977E-09 \)&
 \( 9.45008428970767E-09 \) &
\( 1 \)\\
\hline 
\end{tabular}

\par{} \vspace{0.3cm}\centering

\vspace{0.3cm}\centering 

Bounds using Algorithm II

\begin{tabular}{|c|c|c|c|}
\hline 
Range zeros &
Starting&
NR bounds on Range &
Number of \\
at precision 7&
Value&
zeros at precision 7&
iterations\\
&
&
&
\( NR \) + Algorithm I\\
\hline 
\( 4.993716 \)&
\( 1E-04 \)&
\( 2.09016807585603-05 \) &
\( 4+1 \)\\
\hline 
\( 4.0080655 \)&
\( 1E-04 \)&
\( 3.64417392553444E-05 \)&
\( 4+1 \)\\
\hline 
\( 2.997721 \)&
\( 0.01 \)&
\( 1.6599262404165E-05 \)&
\( 4+1 \)\\
\hline 
\( 1.0 \)&
\( 0.01 \)&
\( 7.33067508160057E-07 \) &
\( 4+1 \)\\
\hline 
\( 0.009422683 \)&
\( 0.01 \)&
\( 1.26764598079487E-08 \)&
\( 4+1 \)\\
\hline 
\( 0.007075149 \)&
\( 1E-10 \)&
 \( 9.45098509041444E-09 \)&
\( 3+1 \)\\
\hline 
\end{tabular}

\par{} \vspace{0.3cm}\centering

\vspace{0.3cm}
{\centering \begin{tabular}{|c|}
\hline 
Zeros at precision 16\\
\hline 
\end{tabular}\par}
\vspace{0.3cm}

\vspace{0.3cm}\centering 

Bounds using Algorithm I for \( \epsilon =0.00000001 \)

\begin{tabular}{|c|c|c|c|}
\hline 
Range zeros &
Value Of &
 Our bounds on Range&
Number of\\
at precision 16&
\( q(0) \)&
 zeros at precision 16&
iterations\\
\hline 
\( 4.99371584412958 \)&
\( 1.79904361618406E-14 \)&
 \( 1.7990436341745E-14 \)&
\( 1 \)\\
\hline 
\( 4.00806551632572 \)&
\( 3.09602233507163E-14 \)&
 \( 3.09602236603186E-14 \)&
\( 1 \)\\
\hline 
\( 2.99772080758747 \)&
\( 1.38592484031255E-14 \)&
 \( 1.3859248541718E-14 \)&
\( 1 \)\\
\hline 
\( 1.0 \)&
\( 5.77846972093955E-16 \)&
 \( 5.77846977872427E-16 \)&
\( 1 \)\\
\hline 
\( 0.00942268285074248 \)&
\( 1.13201138948718E-17 \)&
 \( 1.1320114008073E-17 \)&
\( 1 \)\\
\hline 
\( 0.00707514910648869 \)&
\( 8.61633999170187E-18 \)&
 \( 8.61634007786531E-18 \) &
\( 1 \)\\
\hline 
\end{tabular}

\par{} \vspace{0.3cm}\vspace{0.3cm}\centering 

Bounds using Algorithm II

\begin{tabular}{|c|c|c|c|}
\hline 
Range zeros &
Starting&
 NR bounds on Range&
Number of\\
at precision 16&
value&
 zeros at precision 16&
iterations\\
&
&
&
\( NR \) + Algorithm I\\
\hline 
\( 4.993715844129577 \)&
\( 1E-13 \)&
\( 1.79904363417461E-14 \)&
\( 2+1 \)\\
\hline 
\( 4.008065516325719 \)&
\( 1E-22 \)&
\( 3.09602236603224E-14 \)&
\( 2+1 \)\\
\hline 
\( 2.997720807587473 \)&
\( 1E-09 \)&
\( 1.38592485417189E-14 \)&
\( 2+1 \)\\
\hline 
\( 1.0 \)&
\( 1E-04 \)&
\( 5.77846977872429E-16 \) &
\( 3+1 \)\\
\hline 
\( 0.009422682850742485 \)&
\( 1E-04 \)&
\( 1.13201140080731E-17 \)&
\( 3+1 \)\\
\hline 
\( 0.007075149106488685 \)&
\( 1E-22 \)&
 \( 8.61634007786535E-18 \)&
\( 2+1 \)\\
\hline 
\end{tabular}

\par{} \vspace{0.3cm}\centering

\subsection*{Example 4}

\( g(z)=1000000z^{4}+223069z^{3}-41948404z^{2}+68883845z+125362605 \)

\vspace{0.3cm}
{\centering \begin{tabular}{|c|}
\hline 
Zeros at precision 7\\
\hline 
\end{tabular}\par}
\vspace{0.3cm}

\vspace{0.3cm}\centering

Bounds using Algorithm I for \( \epsilon =0.01 \)

\begin{tabular}{|c|c|c|c|}
\hline 
Range zeros at &
Value at &
 Our bounds on Range&
Number of\\
precision 7&
\( q(0) \)&
zeros at precision 7&
iterations\\
\hline 
\( 4.00102 \)&
\( 1.14317714827099E-03 \)&
 \( 2.9420093174445E-3 \)&
\( 95 \)\\
\hline 
\( 3.998911 \)&
\( 1.1428874099077E-03 \)&
 \( 2.93261103460139E-3 \)&
\( 95 \)\\
\hline 
\( -1.1 \)&
\( 2.11036257663918E-07 \)&
 \( 2.13146620240558E-7 \)&
\( 1 \)\\
\hline 
\( -7.1229995 \)&
\( 5.45032569566617E-07 \)&
 \( 5.50482895262284E-7 \) &
\( 1 \)\\
\hline 
\end{tabular}

\par{} \vspace{0.3cm}\centering

\vspace{0.3cm}\centering

Bounds using Algorithm II

\begin{tabular}{|c|c|c|c|}
\hline 
Range zeros at &
Starting&
 NR bounds on Range&
Number of\\
precision 7&
Value&
zeros at precision 7&
iterations\\
&
&
&
\( NR \) + Algorithm I\\
\hline 
\( 4.00102 \)&
\( 1E-9 \)&
 \( 2.93304493819802E-3 \)&
\( 9+1 \)\\
\hline 
\( 3.998911 \)&
\( 1E-1 \)&
 \( 2.93288220516119E-3 \)&
\( 9+1 \)\\
\hline 
\( -1.1 \)&
\( 1E-5 \)&
 \( 2.11057417059451E-7 \)&
\( 3+1 \)\\
\hline 
\( -7.1229995 \)&
\( 1E-5 \)&
 \( 5.45087277340316E-7 \) &
\( 3+1 \)\\
\hline 
\end{tabular}

\par{} \vspace{0.3cm}\centering

\pagebreak

\vspace{0.3cm}
{\centering \begin{tabular}{|c|}
\hline 
Zeros at precision 16\\
\hline 
\end{tabular}\par}
\vspace{0.3cm}

\vspace{0.3cm} \centering

Bounds using Algorithm I for \( \epsilon =0.00000001 \)

\begin{tabular}{|c|c|c|c|}
\hline 
Range zeros at &
Value at &
 Our bounds on Range&
Number of\\
precision 16&
\( q(0) \)&
zeros at precision 16&
iterations\\
\hline 
\( 4.001019657681204 \)&
\( 1.47961652695472E-12 \)&
 \( 1.47961654175088E-12 \)&
\( 1 \)\\
\hline 
\( 3.998910834038721 \)&
\( 1.47927388424087E-12 \)&
 \( 1.47927389903361E-12 \)&
\( 1 \)\\
\hline 
\( -1.100000011624190 \)&
\( 2.3844189558305E-16 \)&
 \( 2.34844191906748E-16 \)&
\( 1 \)\\
\hline 
\( -7.122999480095735 \)&
\( 6.67094208084199E-16 \)&
 \( 6.67094214755143-16 \) &
\( 1 \)\\
\hline 
\end{tabular}

\par{} \vspace{0.3cm}\centering

\vspace{0.3cm} \centering

Bounds using Algorithm II

\begin{tabular}{|c|c|c|c|}
\hline 
Range zeros at &
Starting&
NR bounds on Range&
Number of\\
precision 16&
Value&
zeros at precision 16&
iterations\\
&
&
&
\( NR \) + Algorithm I\\
\hline 
\( 4.001019657681204 \)&
\( 1E-11 \)&
 \( 1.47961654279053E-12 \)&
\( 3+1 \)\\
\hline 
\( 3.998910834038721 \)&
\( 1E-22 \)&
 \( 1.47927390007277E-12 \)&
\( 2+1 \)\\
\hline 
\( -1.100000011624190 \)&
\( 1E-11 \)&
 \( 2.34844191906748E-16 \)&
\( 2+1 \)\\
\hline 
\( -7.122999480095735 \)&
\( 1E-5 \)&
 \( 6.67094214755144E-16 \) &
\( 3+1 \)\\
\hline 
\end{tabular}

\par{} \vspace{0.3cm}\centering

\subsection*{Example 5:}

\( g(z)=8.7029z^{9}-167z^{8}+463.33z^{7}+1126.1z^{6}+76.241z^{5}-7.0508z^{4}-4085.4z^{3}-1036.1z^{2}-99.729z-54.649 \)

\par{} \vspace{0.7cm}\centering

\vspace{0.3cm}
{\centering \begin{tabular}{|c|}
\hline 
Zeros at precision 7\\
\hline 
\end{tabular}\par}
\vspace{0.3cm}

Bounds using Algorithm I for \( \epsilon =0.0001 \)

\vspace{0.3cm}
{\centering \begin{tabular}{|c|c|c|c|c|}
\hline 
Actual &
Range zeros &
Value of&
Our bounds on &
Number \\
zeros&
at &
\( q(0) \)&
Range zeros &
of\\
&
precision 7&
&
at precision 7&
iterations\\
\hline 
\( 15.0911133 \)&
\( 15.09111 \)&
\( 9.076342589 \)&
\( 9.07725022 \)&
\( 1 \)\\
&
&
\( 63854E-06 \)&
\( 389751E-06 \)&
\\
\hline 
\( 5.60181709 \)&
\( 5.601817 \)&
\( 5.206365229 \)&
\( 5.20688586 \)&
\( 1 \)\\
&
&
\( 73372E-06 \)&
\( 625669E-06 \)&
\\
\hline 
\( 1.471636 \)&
\( 1.471636 \)&
\( 7.629781518 \)&
\( 7.6305444 \)&
\( 1 \)\\
\hline 
&
&
\( 41216E-07 \)&
\( 96564E-07 \)&
\\
\hline 
\( 0.029598+ \)&
\( 0.02959805 \)&
\( 8.494188735 \)&
\( 8.49503815 \)&
\( 1 \)\\
\( i0.204796 \)&
\( +i0.2047964 \)&
\( 02682E-08 \)&
\( 390032E-08 \)&
\\
\hline 
\( 0.029598- \)&
\( 0.02959805 \)&
\( 8.494188735 \)&
\( 8.49503815 \)&
\( 1 \)\\
\( i0.204796 \)&
\( -i0.2047964 \)&
\( 02682E-08 \)&
\( 390032E-08 \)&
\\
\hline 
\( -0.3115666 \)&
\( -0.3115659 \)&
\( 1.57728567 \)&
\( 1.5774434 \)&
\( 1 \)\\
&
&
\( 206155E-07 \)&
\( 0062876E-07 \)&
\\
\hline 
\( -0.396041 \)&
\( -0.396041 \)&
\( 7.26872135 \)&
\( 7.26944822 \)&
\( 1 \)\\
\( +i1.3425333 \)&
\( +i1.342533 \)&
\( 734945E-07 \)&
\( 958519E-07 \)&
\\
\hline 
\( -0.396041 \)&
\( -0.3960412 \)&
\( 7.26872135 \)&
\( 7.26944822 \)&
\( 1 \)\\
\( -i1.3425333 \)&
\( -i1.342533 \)&
\( 734945E-07 \)&
\( 958519E-07 \)&
\\
\hline 
\( -1.9311091 \)&
\( -1.931109 \)&
\( 1.03460887 \)&
\( 1.034712335 \)&
\( 1 \)\\
&
&
\( 418129E-06 \)&
\( 06871E-06 \)&
\\
\hline 
\end{tabular}\par}
\vspace{0.3cm}

\par{} \vspace{0.7cm}\centering

Bounds using Algorithm II

\vspace{0.3cm}
{\centering \begin{tabular}{|c|c|c|c|}
\hline 
Range zeros &
Starting &
Our bounds on &
Number \\
at &
Value&
Range zeros &
of\\
precision 7&
&
at precision 7&
iterations\\
&
&
&
\( NR+ \)Algorithm I\\
\hline 
\( 15.09111 \)&
\( IE-06 \)&
\( 9.07733851 \)&
\( 3+1 \)\\
&
&
\( 508154E-06 \)&
\\
\hline 
\( 5.601817 \)&
\( IE-06 \)&
\( 5.2069566 \)&
\( 3+1 \)\\
&
&
\( 1442314E-06 \)&
\\
\hline 
\( 1.471636 \)&
\( IE-11 \)&
\( 7.63058385 \)&
\( 3+1 \)\\
&
&
\( 852901E-07 \)&
\\
\hline 
\( 0.02959805 \)&
\( IE-11 \)&
\( 8.49504782 \)&
\( 3+1 \)\\
\( +i0.2047964 \)&
&
\( 584932E-08 \)&
\\
\hline 
\( 0.02959805 \)&
\( IE-02 \)&
\( 8.49504782 \)&
\( 4+1 \)\\
\( -i0.2047964 \)&
&
\( 584932E-08 \)&
\\
\hline 
\( -0.3115659 \)&
\( IE-05 \)&
\( 1.57744695 \)&
\( 3+1 \)\\
&
&
\( 490493E-07 \)&
\\
\hline 
\( -0.396041 \)&
\( IE-02 \)&
\( 7.2694861 \)&
\( 4+1 \)\\
\( +i1.342533 \)&
&
\( 2010097E-07 \)&
\\
\hline 
\( -0.3960412 \)&
\( IE-03 \)&
\( 7.2694861 \)&
\( 4+1 \)\\
\( -i1.342533 \)&
&
\( 2010097E-07 \)&
\\
\hline 
\( -1.931109 \)&
\( IE-03 \)&
\( 1.03471867 \)&
\( 4+1 \)\\
&
&
\( 334759E-06 \)&
\\
\hline 
\end{tabular}\par}
\vspace{0.3cm}

\par{} \vspace{0.7cm}\centering

\pagebreak

\vspace{0.3cm}
{\centering \begin{tabular}{|c|}
\hline 
Zeros at precision 16\\
\hline 
\end{tabular}\par}
\vspace{0.3cm}

Bounds using Algorithm I for \( \epsilon =0.00000001 \)

\vspace{0.3cm}
{\centering \begin{tabular}{|c|c|c|c|}
\hline 
Range zeros &
Value of&
Our bounds on &
Number\\
at precision 16&
\( q(0 \))&
Range zeros &
of\\
&
&
at precision 16&
iterations\\
\hline 
\( 15.09111310576714 \)&
\( 5.53750340 \)&
\( 5.53750346 \)&
\( 1 \)\\
&
\( 672293E-15 \)&
\( 209812E-15 \)&
\\
\hline 
\( 5.601817329197456 \)&
\( 3.14720692 \)&
\( 3.14720695 \)&
\( 1 \)\\
&
\( 264054E-15 \)&
\( 411265E-15 \)&
\\
\hline 
\( 1.471636442949784 \)&
\( 5.8533007 \)&
\( 5.85330076 \)&
\( 1 \)\\
&
\( 612049E-16 \)&
\( 12049E-16 \)&
\\
\hline 
\( 0.02959805472609926 \)&
\( 6.97248888 \)&
\( 6.972488954 \)&
\( 1 \)\\
\( +i0.2047964459305099 \)&
\( 506532E-17 \)&
\( 79027E-17 \)&
\\
\hline 
\( 0.02959805472609926 \)&
\( 6.97248888 \)&
\( 6.972488954 \)&
\( 1 \)\\
\( -i0.2047964459305099 \)&
\( 506532E-17 \)&
\( 79027E-17 \)&
\\
\hline 
\( -0.3115658643893903 \)&
\( 1.461254491 \)&
\( 1.46125450 \)&
\( 1 \)\\
&
\( 76875E-16 \)&
\( 63813E-16 \)&
\\
\hline 
\( -0.3960411556175865 \)&
\( 5.65575203 \)&
\( 5.655752091 \)&
\( 1 \)\\
\( +i1.342532717206571 \)&
\( 517524E-16 \)&
\( 73281E-16 \)&
\\
\hline 
\( -0.3960411556175865 \)&
\( 5.65575203 \)&
\( 5.655752091 \)&
\( 1 \)\\
\( -i1.342532717206571 \)&
\( 517524E-16 \)&
\( 73281E-16 \)&
\\
\hline 
\( -1.931108848212619 \)&
\( 7.35055936 \)&
\( 7.350559365 \)&
\( 1 \)\\
&
\( 523337E-16 \)&
\( 23337E-16 \)&
\\
\hline 
\end{tabular}\par}
\vspace{0.3cm}

\par{} \vspace{0.7cm}\centering

Bounds using Algorithm II

\vspace{0.3cm}
{\centering \begin{tabular}{|c|c|c|c|}
\hline 
Range zeros &
Starting&
Our bounds &
Number\\
 at precision 16&
Value&
at precision 16&
of\\
&
&
&
iterations\\
&
&
&
\( NR \) + Algorithm I\\
\hline 
\( 15.091113105767141 \)&
\( 1E-14 \)&
\( 5.53750346209816E-15 \)&
\( 2+1 \)\\
\hline 
\( 5.6018173291974556 \)&
\( 1E-7 \)&
\( 3.14720695411267E-15 \)&
\( 3+1 \)\\
\hline 
\( 1.471636442949784 \)&
\( 0.0001 \)&
\( 5.85330076120492E-16 \)&
\( 3+1 \)\\
\hline 
\( 0.02959805472609926 \)&
\( 1E-22 \)&
\( 6.97248895479027E-17 \)&
\( 2+1 \)\\
\( +i0.2047964459305099 \)&
&
&
\\
\hline 
\( 0.02959805472609926 \)&
\( 1E-14 \)&
\( 6.97248895479027E-17 \)&
\( 2+1 \)\\
\( -i0.2047964459305099 \)&
&
&
\\
\hline 
\( -0.3115658643893903 \)&
\( 1E-7 \)&
\( 1.46125450638131E-16 \)&
\( 3+1 \)\\
\hline 
\( -0.3960411556175865 \)&
\( 0.0001 \)&
\( 5.65575209173284E-16 \)&
\( 3+1 \)\\
\( +i1.342532717206571 \)&
&
&
\\
\hline 
\( -0.3960411556175865 \)&
\( 1E-14 \)&
\( 5.65575209173284E-16 \)&
\( 2+1 \)\\
\( -i1.342532717206571 \)&
&
&
\\
\hline 
\( -1.931108848212619 \)&
\( 1E-22 \)&
\( 7.3505593652334E-16 \)&
\( 2+1 \)\\
\hline 
\end{tabular}\par}
\vspace{0.3cm}

\subsection*{Example 6:}

\( g(z)=z^{20}+10^{12}z^{14}+z^{5}+1 \)

\vspace{0.3cm} \centering 

\begin{tabular}{|c|}
\hline 
Zeros at precision 7  \\
\hline 
\end{tabular}

\par{} \vspace{0.3cm}

Bounds using Algorithm I for \( \epsilon =0.0001 \)

\vspace{0.3cm} \centering 

\begin{tabular}{|c|c|c|c|}
\hline 
Range zeros &
 Value&
 Our bounds on &
 Number\\
 at precision 7&
 of &
 Range zeros&
 of\\
&
 \( q(0) \)&
 at precision 7&
 iterations\\
\hline 
\( 86.60254+i50 \)&
 \( 4.66362417422174E-7 \)&
 \( 4.66410117625052E-7 \)&
 \( 1 \)\\
\hline 
\( 86.60254-i50 \)&
 \( 4.66363481276914E-7 \)&
 \( 4.66410117625041E-7 \)&
 \( 1 \)\\
\hline 
\( 0.1354659+i0.03091968 \)&
 \( 6.1585340861066E-8 \)&
 \( 6.15914993951521E-8 \)&
 \( 1 \)\\
\hline 
\( 0.1354659-i0.03091968 \)&
 \( 6.1585340861066E-8 \)&
 \( 6.15914993951521E-8 \)&
 \( 1 \)\\
\hline 
\( 0.1086348+i0.0866332 \)&
 \( 6.15859858647888E-8 \)&
 \( 6.15921444633753E-8 \)&
 \( 1 \)\\
\hline 
\( 0.1086348-i0.0866332 \)&
 \( 6.15859858647888E-8 \)&
 \( 6.15921444633753E-8 \)&
 \( 1 \)\\
\hline 
\( 0.06028841+i0.1251894 \)&
 \( 6.15851872002481E-8 \)&
 \( 6.15913457189681E-8 \)&
 \( 1 \)\\
\hline 
\( 0.06028841-i0.1251894 \)&
 \( 6.15851872002481E-8 \)&
 \( 6.15913457189681E-8 \)&
 \( 1 \)\\
\hline 
\( 1.666667E-29+i100 \)&
 \( 4.6636348316034E-7 \)&
 \( 4.66410119508656E-7 \)&
 \( 1 \)\\
\hline 
\( 1.666667E-29-i100 \)&
 \( 4.6636348316034E-7 \)&
 \( 4.66410119508656E-7 \)&
 \( 1 \)\\
\hline 
\( -5.140612E-7+i0.1389495 \)&
 \( 6.15855211502818E-8 \)&
 \( 6.15916797023968E-8 \)&
 \( 1 \)\\
\hline 
\( -5.140612E-7-i0.1389495 \)&
 \( 6.15855211502818E-8 \)&
 \( 6.15916797023968E-8 \)&
 \( 1 \)\\
\hline 
\( -0.0602879+i0.125189 \)&
 \( 6.15858807166052E-8 \)&
 \( 6.15920393046768E-8 \)&
 \( 1 \)\\
\hline 
\( -0.0602879-i0.125189 \)&
 \( 6.15858807166052E-8 \)&
 \( 6.15920393046768E-8 \)&
 \( 1 \)\\
\hline 
\( -0.1086354+i0.0866340 \)&
 \( 6.15850868338623E-8 \)&
 \( 6.15912453425457E-8 \)&
 \( 1 \)\\
\hline 
\( -0.1086354-i0.0866340 \)&
 \( 6.15850868338623E-8 \)&
 \( 6.15912453425456E-8 \)&
 \( 1 \)\\
\hline 
\( -0.1354657+i0.03091868 \)&
 \( 6.15857290119851E-8 \)&
 \( 6.15918875848863E-8 \)&
 \( 1 \)\\
\hline 
\( -0.1354657-i0.03091868 \)&
 \( 6.15857290119851E-8 \)&
 \( 6.15918875848863E-8 \)&
 \( 1 \)\\
\hline 
\( -86.60254+i50 \)&
 \( 4.66363481988736E-7 \)&
 \( 4.66410118336934E-7 \)&
 \( 1 \)\\
\hline 
\( -86.60254-i50 \)&
 \( 4.66363481988736E-7 \)&
 \( 4.66410118336934E-7 \)&
 \( 1 \) \\
\hline 
\end{tabular}

\par{} \vspace{0.3cm}

\pagebreak

Bounds using Algorithm II

\vspace{0.3cm} \centering 

\begin{tabular}{|c|c|c|c|}
\hline 
Range zeros &
 Starting &
 NR bounds on &
 Number\\
 at precision 7&
 Value &
 Range zeros&
 of\\
&
&
 at precision 7&
 iterations\\
&
&
&
 \( NR+ \)Algorithm I\\
\hline 
\( 86.60254+i50 \)&
 \( 1E-10 \)&
 \( 4.66410190829221E-7 \)&
 \( 2+1 \)\\
\hline 
\( 86.60254-i50 \)&
 \( 1E-06 \)&
 \( 4.66410190829219E-7 \)&
 \( 2+1 \)\\
\hline 
\( 0.1354659+i0.03091968 \)&
 \( 1E-05 \)&
 \( 6.15920100912469E-8 \)&
 \( 3+1 \)\\
\hline 
\( 0.1354659-i0.03091968 \)&
 \( 1E-06 \)&
 \( 6.15920095777814E-8 \)&
 \( 2+ \)\( 1 \)\\
\hline 
\( 0.1086348+i0.086633225 \)&
 \( 1E-10 \)&
 \( 6.15926551678678E-8 \)&
 \( 2+1 \)\\
\hline 
\( 0.1086348-i0.086633225 \)&
 \( 1E-05 \)&
 \( 6.15926551700748E-8 \)&
 \( 3+1 \)\\
\hline 
\( 0.06028841+i0.1251894 \)&
 \( 1E-05 \)&
 \( 6.15918564125151E-8 \)&
\( 3+1 \)\\
\hline 
\( 0.06028841-i0.1251894 \)&
 \( 1E-03 \)&
 \( 6.15918564125148E-8 \)&
 \( 3+1 \)\\
\hline 
\( 1.666667E-29+i100 \)&
 \( 1E-04 \)&
 \( 4.66410192449928E-7 \)&
 \( 2+1 \)\\
\hline 
\( 1.666667E-29-i100 \)&
 \( 1E-03 \)&
 \( 4.66410192712841E-7 \)&
 \( 3+1 \)\\
\hline 
\( -5.140612E-7+i0.13894955 \)&
 \( 1E-04 \)&
 \( 6.15921904014471E-8 \)&
 \( 3+1 \)\\
\hline 
\( -5.140612E-7-i0.13894955 \)&
 \( 1E-06 \)&
 \( 6.159218988798E-8 \)&
 \( 2+1 \)\\
\hline 
\( -0.0602879+i0.125189 \)&
 \( 1E-03 \)&
 \( 6.15925500096416E-8 \)&
 \( 3+1 \)\\
\hline 
\( -0.0602879-i0.125189 \)&
 \( 1E-06 \)&
 \( 6.15925494961725E-8 \)&
 \( 2+1 \)\\
\hline 
\( -0.1086354+i0.08663403 \)&
 \( 1E-04 \)&
 \( 6.15917560344614E-8 \)&
 \( 3+1 \)\\
\hline 
\( -0.1086354+i0.08663403 \)&
 \( 1E-10 \)&
 \( 6.15917560322547E-8 \)&
 \( 2+1 \)\\
\hline 
\( -0.1354657+i0.03091868 \)&
 \( 1E-06 \)&
 \( 6.15923977739059E-8 \)&
 \( 2+1 \)\\
\hline 
\( -0.1354657-i0.03091868 \)&
 \( 1E-10 \)&
 \( [6.15923982851674E-8 \)&
 \( 2+1 \)\\
\hline 
\( -86.60254+i50 \)&
 \( 1E-04 \)&
 \( 4.66410191278206E-7 \)&
 \( 2+1 \)\\
\hline 
\( -86.60254-i50 \)&
 \( 1E-10 \)&
 \( 4.66410191541114E-7 \)&
 \( 2+1 \) \\
\hline 
\end{tabular}

\par{} \vspace{0.3cm}

\vspace{0.3cm} \centering 

\begin{tabular}{|c|}
\hline 
Zeros at precision 16  \\
\hline 
\end{tabular}

\par{} \vspace{0.3cm}

Bounds using Algorithm I for \( \epsilon =0.00000001 \)

\vspace{0.3cm} \centering 

\begin{tabular}{|c|c|c|c|}
\hline 
Range zeros &
 Value of&
 Our bounds &
 Number\\
 at precision 16&
 \( q(0) \)&
 on Range zeros&
 of\\
&
&
 at precision 16&
 iterations\\
\hline 
\( 86.60254037844386+i50 \)&
 \( 2.77298663308076E-15 \)&
 \( 2.77298666081066E-15 \)&
 \( 1 \)\\
\hline 
\( 86.60254037844386-i50 \)&
 \( 2.77298663308073E-15 \)&
 \( 2.77298666081066E-15 \)&
 \( 1 \)\\
\hline 
\( 0.135465908523133 \)&
 \( 8.26333010777928E-17 \)&
 \( 8.26333019041465E-17 \)&
 \( 1 \)\\
 \( +i0.03091968468442907 \)&
&
&
\\
\hline 
\( 0.135465908523133 \)&
 \( 8.26333010778121E-17 \)&
 \( 8.26333019041465E-17 \)&
 \( 1 \)\\
 \( -i0.03091968468442907 \)&
&
&
\\
\hline 
\( 0.1086348117116894 \)&
 \( 8.26340929550181E-17 \)&
 \( 8.26340937813606E-17 \)&
 \( 1 \)\\
 \( +i0.0866332513523468 \)&
&
&
\\
\hline 
\( 0.1086348117116894 \)&
 \( 8.26340929550181E-17 \)&
 \( 8.26340937813606E-17 \)&
 \( 1 \)\\
 \( -i0.0866332513523468 \)&
&
&
\\
\hline 
\( 0.06028841320891334 \)&
 \( 8.26331054824203E-17 \)&
 \( 8.26331063087529E-17 \)&
 \( 1 \)\\
 \( +i0.125189441295516 \)&
&
&
\\
\hline 
\( 0.06028841320891334 \)&
 \( 8.26331054824203E-17 \)&
 \( 8.26331063087529E-17 \)&
 \( 1 \)\\
 \( -i0.125189441295516 \)&
&
&
\\
\hline 
\( 1.666666666666667E-29+i100 \)&
 \( 2.77298663308073E-15 \)&
 \( 2.77298666081066E-15 \)&
 \( 1 \)\\
\hline 
\( 1.66666666666667E-29-i100 \)&
 \( 2.77298663308073E-15 \)&
 \( 2.77298666081066E-15 \)&
 \( 1 \)\\
\hline 
\( -5.140611950289672E-7 \)&
 \( 8.26335449713578E-17 \)&
 \( 8.26335457976949E-17 \)&
 \( 1 \)\\
 \( +i0.1389495494401665 \)&
&
&
\\
\hline 
\( -5.14061950289672E-7 \)&
 \( 8.26335449713578E-17 \)&
 \( 8.26335457976949E-17 \)&
 \( 1 \)\\
 \( -i0.1389495494401665 \)&
&
&
\\
\hline 
\( -0.06028748690264759 \)&
 \( 8.26339844245486E-17 \)\( 1 \)&
 \( 8.26339852508899E-17 \)&
 \( 1 \)\\
 \( +i0.1251889952099292 \)&
&
&
\\
\hline 
\( -0.06028748690264759 \)&
 \( 8.26339844245486E-17 \)&
 \( 8.26339852508899E-17 \)&
 \( 1 \)\\
 \( -i0.1251889952099292 \)&
&
&
\\
\hline 
\( -0.1086354527355147 \)&
 \( 8.26329969321128E-17 \)&
 \( 8.26329977584442E-17 \)&
 \( 1 \)\\
 \( +i0.08663402895368696 \)&
&
&
\\
\hline 
\( -0.1086354527355147 \)&
 \( 8.26329969321128E-17 \)&
 \( 8.26329977584442E-17 \)&
 \( 1 \)\\
 \( -i0.08663402895368696 \)&
&
&
\\
\hline 
\( -0.135465797443785 \)&
 \( 8.26337888538945E-17 \)&
 \( 8.2633789680234E-17 \)&
 \( 1 \)\\
 \( +i0.03091868233921389 \)&
&
&
\\
\hline 
\( -0.1354657443785 \)&
 \( 8.26337888538945E-17 \)&
 \( [8.2633789680234E-17 \)&
 \( 1 \)\\
 \( -i0.03091868233921389 \)&
&
&
\\
\hline 
\( -86.60254037844386+i50 \)&
 \( 2.77298663308073E-15 \)&
 \( 2.77298666081066E-15 \)&
 \( 1 \)\\
\hline 
\( -86.60254037844386+i50 \)&
 \( 2.77298663308073E-15 \)&
 \( 2.77298666081066E-15 \)&
 \( 1 \) \\
\hline 
\end{tabular}

\par{} \vspace{0.3cm}

\pagebreak

Bounds using Algorithm II

\vspace{0.3cm} \centering 

\begin{tabular}{|c|c|c|c|}
\hline 
Range zeros&
 Starting&
 NR bounds &
 Number\\
 at precision 16&
 Value&
 on Range zeros&
 of \\
&
&
 at precision 16&
 iterations\\
&
&
&
 \( NR+ \)Algorithm I\\
\hline 
\( 86.60254037844386+i50 \)&
 \( 1E-15 \)&
 \( 2.77298666081066E-15 \)&
 \( 2+1 \)\\
\hline 
\( 86.60254037844386-i50 \)&
 \( 1E-10 \)&
 \( 2.77298666081066E-15 \)&
 \( 2+1 \)\\
\hline 
\( 0.135465908523133 \)&
 \( 1E-07 \)&
 \( 8.26333018961745E-17 \)&
 \( 2+1 \)\\
 \( +i0.03091968468442907 \)&
&
&
\\
\hline 
\( 0.135465908523133 \)&
 \( 1E-15 \)&
 \( 8.26333019041475E-17 \)&
 \( 1+1 \)\\
 \( -i0.03091968468442907 \)&
&
&
\\
\hline 
\( 0.1086348117116894 \)&
 \( 1E-10 \)&
 \( 8.26340937813615E-17 \)&
 \( 2+1 \)\\
 \( +i0.0866332513523468 \)&
&
&
\\
\hline 
\( 0.1086348117116894 \)&
 \( 1E-16 \)&
 \( 8.26340937813615E-17 \)&
 \( 1+1 \)\\
 \( -i0.0866332513523468 \)&
&
&
\\
\hline 
\( 0.06028841320891334 \)&
 \( 1E-10 \)&
 \( 8.26331063087538E-17 \)&
 \( 2+1 \)\\
 \( +i0.125189441295516 \)&
&
&
\\
\hline 
\( 0.06028841320891334 \)&
 \( 1E-07 \)&
 \( 8.26331063007808E-17 \)&
 \( 2+1 \)\\
 \( -i0.125189441295516 \)&
&
&
\\
\hline 
\( 1.666666666666667E-29+i100 \)&
 \( 1E-15 \)&
 \( 2.77298666081066E-15 \)&
 \( 2+1 \)\\
\hline 
\( 1.666666666666667E-29-i100 \)&
 \( 1E-10 \)&
 \( 2.77298666081066E-15 \)&
 \( 2+1 \)\\
\hline 
\( -5.140611950289672E-7 \)&
 \( 1E-10 \)&
 \( 8.26335457976958E-17 \)&
 \( 2+1 \)\\
 \( +i0.1389495494401665 \)&
&
&
\\
\hline 
\( -5.14061950289672E-7 \)&
 \( 1E-07 \)&
 \( 8.26335457897226E-17 \)&
 \( 2+1 \)\\
\hline 
\( -i0.1389495494401665 \)&
&
&
\\
\hline 
\( -0.06028748690264759 \)&
 \( 1E-16 \)&
 \( 8.26339852508908E-17 \)&
 \( 1+1 \)\\
 \( +i0.1251889952099292 \)&
&
&
\\
\hline 
\( -0.06028748690264759 \)&
 \( 1E-07 \)&
 \( 8.26339852429179E-17 \)&
 \( 2+1 \)\\
 \( -i0.1251889952099292 \)&
&
&
\\
\hline 
\( -0.1086354527355147 \)&
 \( 1E-15 \)&
 \( 8.26329977584451E-17 \)&
 \( 1+1 \)\\
\hline 
\( +i0.08663402895368696 \)&
&
&
\\
\hline 
\( -0.1086354527355147 \)&
 \( 1E-16 \)&
 \( 8.26329977584451E-17 \)&
 \( 1+1 \)\\
 \( -i0.08663402895368696 \)&
&
&
\\
\hline 
\( -0.135465797443785 \)&
 \( 1E-10 \)&
 \( 8.26337896802349E-17 \)&
 \( 2+1 \)\\
\hline 
\( +i0.03091868233921389 \)&
&
&
\\
\hline 
\( -0.135465797443785 \)&
 \( 1E-15 \)&
 \( 8.26337896802349E-17 \)&
 \( 1+1 \)\\
\hline 
\( -i0.03091868233921389 \)&
&
&
\\
\hline 
\( -86.60254037844386+i50 \)&
 \( 1E-16 \)&
 \( 2.77298666081066E-15 \)&
 \( 2+1 \)\\
\hline 
\( -86.60254037844386-i50 \)&
 \( 1E-07 \)&
 \( 2.77298666081066E-15 \)&
 \( 2+1 \) \\
\hline 
\end{tabular}

\par{} \vspace{0.3cm}

}
% EXPLANATION OF RESULTS ISOLATED ZEROS EX 1 and 2

\section{Results and observations}

We tabulate the bounds computed for examples 1 through 6 
using our Algorithms I and II
based on Theorem \ref{maintheorem}. We also 
tabulate the number of iterations
required by Algorithms I and II. For Algorithm II, the Newton-Raphson 
iterations are followed by steps of Algorithm I. So, the numbers 
of iterations of both these stages are depicted in the tables.
Note that the value of $q(0)$ for each approximate zero is
very close to (and slightly lesser than) 
the error bound, as argued in Section \ref{algoI}.
 
For examples 1 and 2 the actual zeros are known.
Both these polynomials (see \cite{smith}) 
have distinct and well separated zeros.
Approximate zeros as computed by the ZERPOL algorithm \cite{smith0} are
used in \cite{smith}; Smith computes error bounds for these zeros using his 
method (as in \cite{smith}) based on Gerschgorin's theorems.
Error bounds computed by our Algorithm I using Theorem \ref{maintheorem}
for these ZERPOL zeros are 
comparatively inferior to those of \cite{smith}. For example 1,
our bounds are comparable to those of Smith; the bounds on 
errors in the last two zeros are very close. For example 2, 
our bounds  for the seventh and eighth zeros are in fact better than 
those obtained in \cite{smith}. The bounds for the other eight zeros are almost 
of the same order in our case whereas Smith's bound for these eight zeros are
much better. However, 
our method is scalable. We show that our bounds are sharper 
for more accurate approximations. As we increase the number of significant 
digits in the approximate zeros computed using the {\it Range} 
software \cite{aberth}, we find that our error bounds too improve as expected.
We use the {\it Range} software \cite{aberth} 
for computing approximate zeros for all examples 
at precisions of 7 and 16 decimal significant digits, respectively. 

Our method works well for polynomials with close zeros as shown in 
examples 3 and 4. The close zeros in example 3 are smaller than the 
rest of the zeros. Error bounds computed using our Algorithms I and II 
for zeros computed using Range \cite{aberth} at the two 
precision values of 7 and 16 
are presented. The bounds in example 3 for small and 
close zeros are much
sharper than the bounds for the 
other zeros because the numerator $l(r)$ of $q(r)$
becomes smaller for smaller zeros. In contrast, we see a coarser 
bound for close zeros with larger magnitudes in example 4, even though 
they are as close as the close zeros in example 3. 

The bounds for the two close zeros (first two zeros)
in example 4 at precision 7 
are larger than the separation between the zeros. 
The two circles with radii equal to 
these two computed bounds and centred at the respective approximate 
zeros enclose both the approximate zeros. 
This is guaranteed by Rouche's theorem because the condition
$r>q(r)$ holds for the bound $r$ for these two zeros. In other words,
the two approximate zeros and the two corresponding actual zeros lie in
both the above mentioned circles. 
Our bounds for the more separated
zeros in example 4 are better that those for the closer zeros.
Increasing the accuracy of approximation helps in getting sharp bounds
isolating each zero. This is observed even for examples 3 and 4 
at precision 16. 

Example 5 is from \cite{henrici} and example 6 is from \cite{posso}.
These are relatively higher degree polynomials and our Algorithms I and II
show consistent and good bounds. We make an crucial observation over all 
examples; we note that the relative error in approximate zeros as computed by
our Algorithms I and II are roughly of the same order of magnitude over all
zeros for each example.

\section{Conclusions}

Posteriori error bounds on approximate zeros of univariate polynomials may be 
used in geometric computations with high degree curves and surfaces
as well as in 
various scientific, engineering and computational mathematics applications.
Our method for computing error bounds can be used in 
algorithms for computing approximate zeros for polynomials to 
desired accuracies; for instance, a Newton-Raphson based algorithm for 
computing approximate zeros may repeatedly compute error bounds in  
each iteration using our 
method to decide whether approximate zeros of desired accuracy have 
already been computed.
Use of high precision 
in zeros' computations can help generating as close approximations to actual 
zeros as one wishes. Krishnan et al. \cite{manocha2001} use Durand-Kerner
iterations for fast convergence to all zeros of a polynomial simultaneously; 
at each step they check the quality of approximation achieved up to that step
using Smith's method \cite{smith} for computing upper bounds on errors in the 
approximations of zeros obtained. Our method of computing similar error
bounds as shown in this paper can be used in the place of Smith's method.

We have used requisite high precision in our algorithms 
for decision making steps involving 
inequalities and in the evaluation of polynomials of high degree.
Our method works for polynomials with 
close clustered zeros; the problem in such cases is that we might require
larger numbers of significant digits in the approximate zeros for computing 
good bounds. With coarser approximations, we might get poorer bounds (see 
example 4). However, for the case of close roots with small magnitudes, as in
example 3, we may get sharp bounds even with coarse approximations of 
zeros. We can also extend our method for handling the case of polynomials 
whose zeros have multiplicity greater than unity. The 
bounds obtained for such cases are likely to be coarse and we feel that only 
very precise approximations of zeros would give sharper bounds.


\begin{thebibliography}{10}
\bibitem{aberth} Aberth O. and Schaefer M. J., ``Range arithmetic
software. {\it http://www.math.tamu.edu/~oliver.aberth/soft.html}''.
\bibitem{ahlfors}Ahlfors, L. V., ``Complex Analysis", McGraw-Hill International Edition, 1979.
\bibitem{exactleda} Burnikel C.,  Könemann J., Mehlhorn K., 
Näher S., Schirra S., and Uhrig C., 
"Exact geometric computation in LEDA", Proc. Symp. on Compu. Geom. 1995,
pages C18-C19. 
\bibitem{henricivol1} Henrici P., ``Applied and Computational Complex 
Analysis, Vol. 1'', John Wiley and Sons, New York, 1974. 
\bibitem{henrici}Henrici, P., AND Watkins, B.O., ``Finding zeros of 
a polynomials by the Q-D
algorithm'', Comm. ACM 8,9 (Sept 1965), 570-574.
\bibitem{JT72} Jenkins M. A. and Traub J. F., 
Algorithm 419: Zeros of a complex polynomial, Comm. ACM 15 (1972), 97-99.
\bibitem{koul}Koul, Rakesh, ``A System for the Exact Computation of Orientation of Transformed
Geometric Objects'', Masters Thesis submitted to the Department of Computer
Science and Engineering, Indian Institute of Technology, Kharagpur, 2000.
\bibitem{manocha2001} Krishnan S., Foskey M., Culver T., Keyser J and
Manocha D., "PRECISE: Efficient multiprecision evaluation of 
algebraic roots and predicates for reliable geometric 
computation", Proc. ACM Symp. on Comput. Geom., June 2001, 274-283.
\bibitem{M73} Madsen K., 
A root-finding algorithm based on Newton's method, BIT 13 (1973), 71-75.
\bibitem{mehlhornnaher99} Mehlhorn K. and Naher S., 
``The LEDA platform for Combinatorial 
and Geometric Computing'', Cambridge University Press, 1999.
\bibitem{mpvt96} Mukherjee, M., Pal, S. P., Varvani M. K.,
Tripathi M., 
``Safe implementation of set operators using 
finite precision'', CSG 96, Proceedings of the conference held in 
Winchester, UK, Set-theoretic Solid Modelling: 
Techniques and Applications, pp. 291-305, April 1996.
\bibitem{posso}PoSSo- Polynomial System Solving project, 
http://www-sop.inria.fr/saga/POL.
\bibitem{pan1}Pan, V.Y., ``Approximating Complex Polynomial Zeros: Modified Weyl's Quad-tree
Construction and Improved Newton's Iteration'', J. of Complexity 16, 213-264
(2000).
\bibitem{pan2}Pan, V.Y., ``Solving a polynomial equation: Some history and recent progress'',
SIAM Review, vol.39, No.2, pp. 187-220, June 1997.
\bibitem{rama}P. H. D. Ramakrishna, ``Bounding errors in the computation of
trigonometric functions and roots of 
polynomials'', M. Tech. (Master's Thesis), Department
of Computer Science and Engineering, Indian Institute of Technology, Kharagpur,
721302, India. 
July, 2000.
\bibitem{techrep} P. H. D. Ramakrishna, S. Bhalla, H. Basu and S. P. Pal,
``Computing bounds on error in solutions of algebraic equations using Rouche's
theorem'', Technical Report \# TR/IIT/CSE/2002/SPP1, Department of Computer Science
and Engineering, Indian Institute of Technology, Kharagpur, 721302, India.
%\bibitem{complex1}Remmert, ``Theory of Complex Functions'', Springer-Verlag, 1991.
\bibitem{smale}Smale, S., ``The Fundamental Theorem of Algebra and Complexity Theory'', Bull.
Amer. Math. Soc., 4 (1981), pp. 1-36.
\bibitem{smith0}Smith, B.T., ``ZERPOL, a zero finding algorithm for 
polynomials using Laguerre's method'', Proc. 1967 Army Numerical 
Analysis Conference, Madison, Wis., 
May 1967 (Rep. 67-3, US Army Res. 
Office- Durham, Durham, N.C., Nov 1967). pp. 153-174.
\bibitem{smith}Smith, B.T., ``Error Bounds for Zeros of a Polynomial Based Upon Gerschgorin's
Theorem'', Journal of ACM, Vol.17, No.4, Oct'1970, pp. 661-674.
\end{thebibliography}
\end{document}